\documentclass[11pt]{amsart}
\newcommand{\N}{\ensuremath{ \mathbf N }}

\usepackage{url}
\usepackage{graphicx}
\usepackage{amssymb}
\usepackage{epstopdf}
\title{Perfect Powers of five with few ternary digits}
\author{Satyanand Singh}
\address{New York City College Of technology of CUNY}
\email{SSingh@citytech.cuny.edu}
\date{\today}                                           
\begin{document}
\begin{abstract}
In this note we will analyze a diophantine equation raised by Michael Bennett in \cite{MB} that is pivotal in establishing that powers of five has few digits in its ternary expansion. We will show that the Diophantine equation $3^{a}+3^{b}+2=n^5$, where $(n,3)=1$ and $a>b>0$ is insoluble for pairs of positive integers $(a,b)$ where they are both even or one is even and the other is odd. In the case where both $(a,b)$ are odd, there is one known solution $2^5=3^3+3^1+2.$ We will show that there are no other solutions to the diophantine equation for $n^{5}<32\left(1+3(10^6)\right)^5$.

\end{abstract}
\maketitle

\section{Introduction}
In \cite{MB}, Bennett showed that the occurrence of at most three digits in the ternary expansion of positive integers raised to certain primes are rare under the conditions stated in the theorem below.\\ 

{\bf Theorem B}. {\it {For any positive integer $n$, with $(n,3)=1,$ if $n^{q}$ has at most three ternary digits for $q$ a prime number when $q=3$ or $7\leq q <1000$, then $n^{q}=13^3=1+3^2+3^7.$}} \\

In classifying the ternary expansions for $q=5$, Bennett dispensed of the equations:
$2^{\delta_{1}}3^{a}+2^{\delta_{2}}=n^5$, for $a>0$, $\delta_{i}\in\{0,1\}$
 and $2^{\delta_{1}}3^{a}+2^{\delta_{2}}3^{b}+2^{\delta_{3}}=n^5$, for $a>b>0$, $\delta_{i}\in\{0,1\}$, except when $(\delta_{1},\delta_{2}, \delta_{3})\neq(0,0,1)$ by showing their insolubility.\\
 
We will establish for $(\delta_{1},\delta_{2}, \delta_{3})=(0,0,1)$, i.e. in the instance when $3^{a}+3^{b}+2=n^5$, that this equation is insoluble when the pairs $(a,b)$ are both even, or one is even and the other is odd. We will use obstructions in the ring $\mathbf{Z}/8\mathbf{Z}$ to establish this claim. We will also use divisibility properties and a maple program to establish that other than the exception, $2^5=3^3+3^1+2$, there are no other powers of five up to $32\left(1+3(10^6)\right)^5$ that are relatively prime to $3$ with at most three digits in its ternary expansion. Note that $\N_{0}$ denote the non-negative integers and $(a, b)=1$ indicates that $a$ and $b$ are relatively prime.

\section{Powers of five and  ternary expansions of the form $3^{a}+3^{b}+2$}
{\bf Lemma 1.}
{\it For $a$ and $b$ even and $a>b\geq2$, the equation $3^{a}+3^{b}+2=n^5$ is insoluble.}
\begin{proof}

Let $a=2t_{1}$, and $b=2t_{2}$, with $t_{1}, t_{2}\in \bf{N_{0}}$ and $t_{1}>t_{2}\geq1.$ It follows that $3^{2t_{1}}+3^{2t_{2}}+2 \equiv \{4\} \mod 8.$ For even $n$, clearly $n^5 \equiv \{0\} \mod 8.$ For $n$ odd, say $n=2r+1$ and $r\in\bf{N_{0}},$ we have that $(2r+1)^5\equiv 2r+1 \mod 8\equiv\{1,3,5,7\} \mod 8.$ In this instance we have an obstruction since $\{4\}\cap \{0,1,3,5,7\}=\emptyset$ and our lemma is proven. 

\end{proof}\

We will now consider the instance when $a$ takes even values and $b$ takes on odd values.\\

{\bf Lemma 2}.
{\it For $a$ even and $b$ odd and $a>b\geq1$, the equation $3^{a}+3^{b}+2=n^5$ is insoluble.}

\begin{proof}
Let $a=2s_{1}$, and $b=2s_{2}+1$, with $s_{1}, s_{2}\in \bf{N_{0}}$ and $s_{1}>s_{2}\geq0.$ It follows that $3^{2s_{1}}+3^{2s_{2}+1}+2\equiv \{6\} \mod 8.$ We saw above that $n^5\equiv\{0,1,3,5,7\} \mod 8.$ We also have an obstruction since $\{6\}\cap \{0,1,3,5,7\}=\emptyset$ and our lemma is proven.
\end{proof}

{\bf Corollary 1}. 
{\it For $a$ odd and $b$ even and $a>b\geq2$, the equation $3^{a}+3^{b}+2=n^5$ is insoluble.} Corollary 1 follows easily from the proof of lemma 2 by a symmetric argument.\\

We will take a few more stabs at using moduli to sieve out some cases for $n$. For the equation $3^{a}+3^{b}+2=n^5$, apply modulo $3$ to both sides. It follows that $2\equiv n^5\mod 3\equiv n\mod 3.$ The latter expression follows from Fermat's little theorem and the stipulation that $(n,3)=1.$ From this we see that $n=3t+2$, where $t\in \bf{N_{0}}.$ Now apply modulo $2$ to both sides of our updated equation, $3^{a}+3^{b}+2=(3t+2)^5$ to get $t^5\mod 2\equiv 0\mod 2.$ It now follows that $t$ is even. We will let $t=2j$, where $j\geq1$ and $j\in\N_{0}$ to obtain $n=6j+2.$
We will now work with $n$ as our starting point, and make use of the Maple software to show that fifth powers of $(6j+2)^5$ has a ternary expansion that consists of at least four non-zero digits. Define $\gamma_{3}(n^5)$ as the number of non-zero digits in the ternary expansion of $n^5$.\\

{\bf Theorem 1}.
The  Diophantine equation $3^{a}+3^{b}+2=n^5$, where $(n,3)=1$ and $a>b>0$, is insoluble for $2<n\leq 2+6(10^6),$ where $n\in\N_{0}.$ \\

\begin{proof}
We implemented code in Maple that computed and plotted the sequence $\gamma_{3}((6j+2)^5)$ for $j=1, 2,...,10^6.$ Our program also plotted on the same axes the horizontal line $(j,3)$ and the two sequences: $\left(j,{\frac{\ln(6j+2)^5}{3\ln3}}\right)$ and $\left(j,{\frac{\ln(6j+2)^5}{\ln3}}\right)$  for this range of $j$ values. We provide by way of example, two plots below for the intervals:$[1, 100000]$ and $[900001, 1000000].$ Our graphs illustrates that $\gamma_{3}((6j+2)^5)>3$ which shows that $3^{a}+3^{b}+2=n^5$ cannot occur for powers of five up to $32\left(1+3(10^6)\right)^5.$ 
\end{proof}

\begin{figure} [ht]
\centering
\begin{tabular}{cc}
\begin{minipage}{2in}
\end{minipage}
\includegraphics[scale=.4]{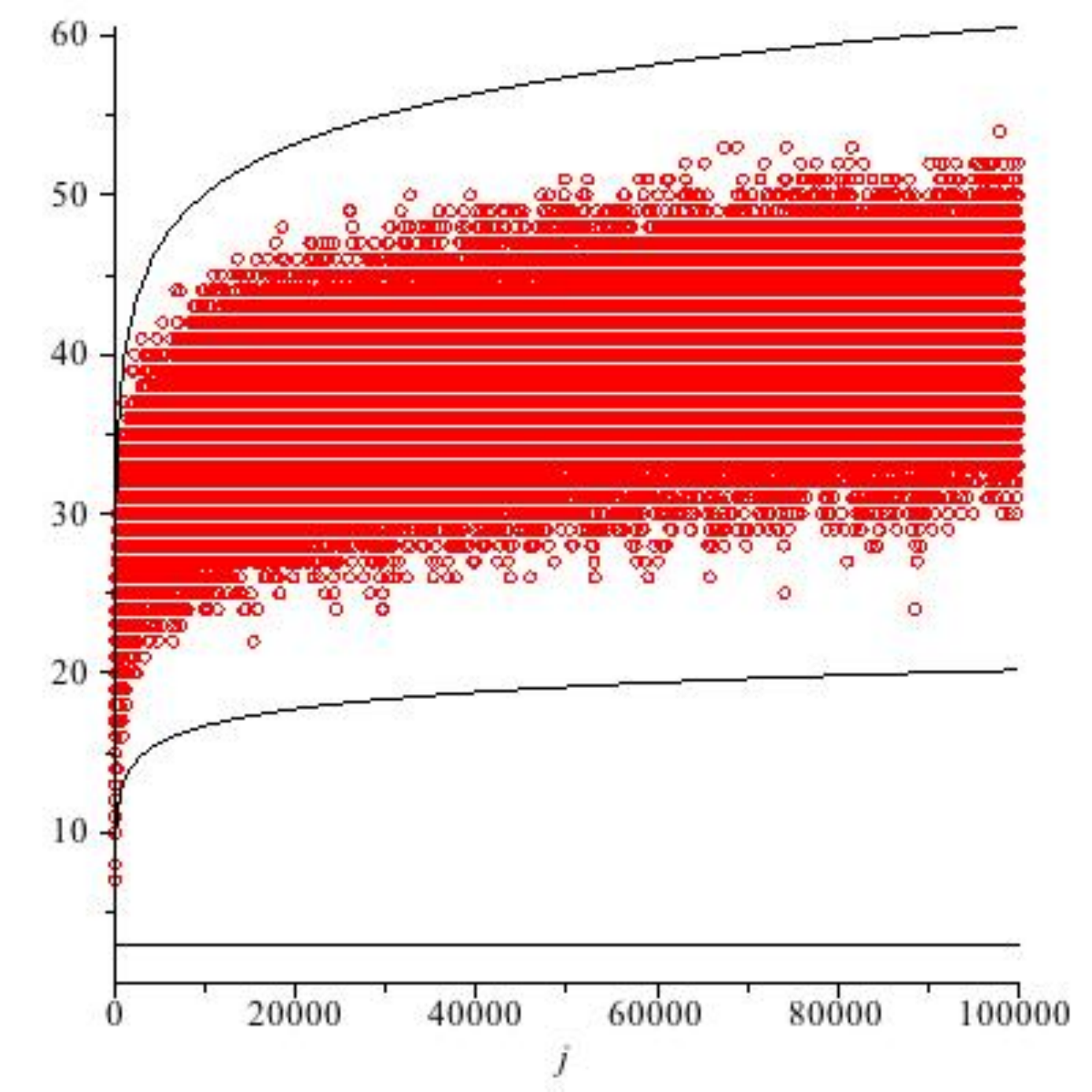}
\includegraphics[scale=.4]{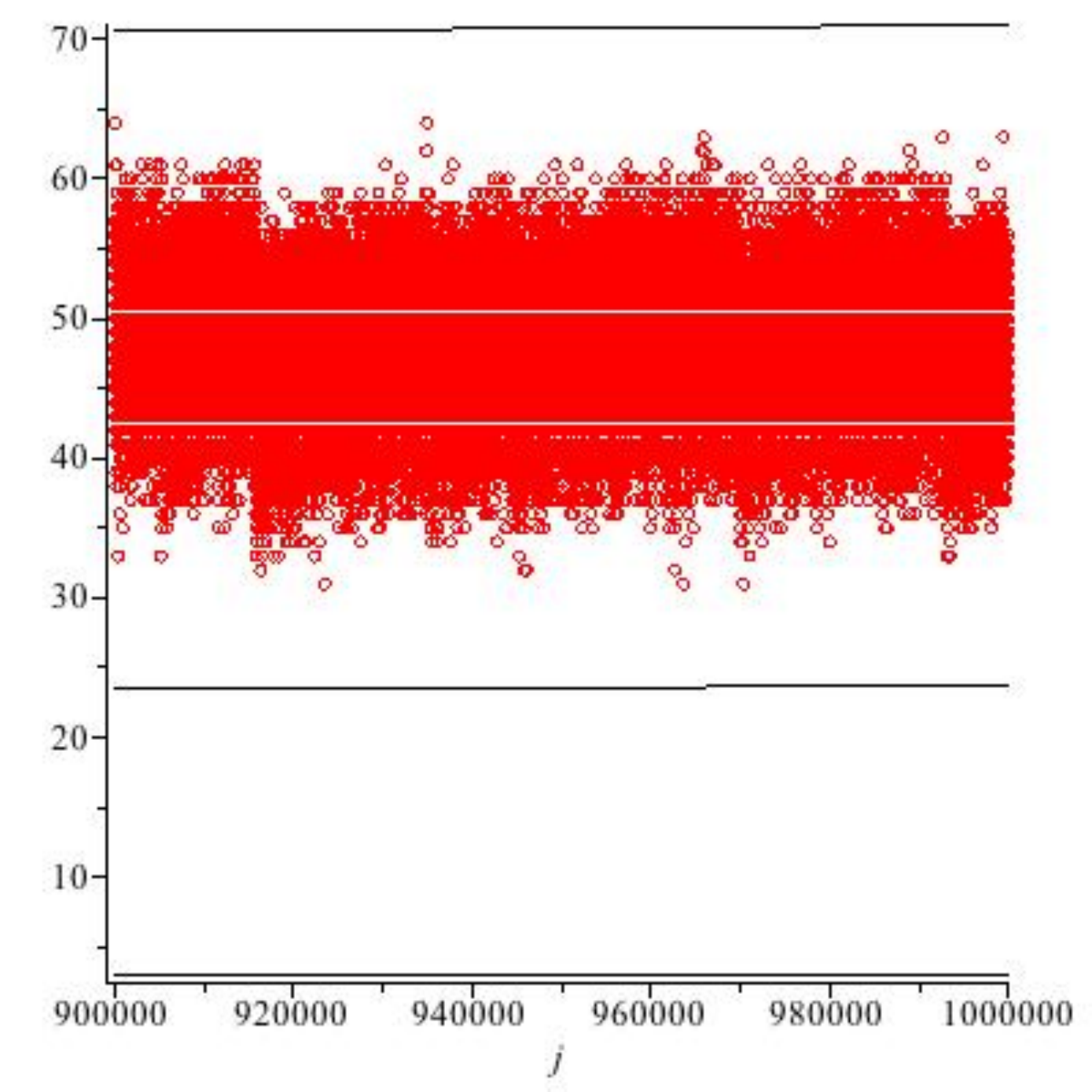}
\end{tabular}
\end{figure}

The two outer graphs which form an envelope around $\gamma_{3}((6j+2)^5)$ illustrate that: 
\[
{\left(\frac{\ln(6j+2)^5}{3\ln3}\right)} <\gamma_{3}((6j+2)^5)<{\left(\frac{\ln(6j+2)^5}{\ln3}\right)}.
\]

The upper bound is easily seen by finding the power of $3$ that is closest to $(6j+2)^5 $ but does not exceed it. We can also say for certain that $\gamma_{3}((6j+2)^5)\geq3,$ since Bennett dispensed of the two term case in \cite{MB} and equality occurs when $2^5=3^3+3^1+2.$ We were not able to prove the lower bound suggested by the experimental results, i.e. $ \gamma_{3}((6j+2)^5)>{\left(\frac{\ln(6j+2)^5}{3\ln3}\right)}$ for $j\geq{1}$. This would completely resolve the case for $q=5.$  

\section{Open Problems and Comments}

We did not resolve the case when both $a$ and $b$ are odd in Bennett's diophantine equation. This case is stated in problem 1 below, and an answer would resolve the $q=5$ case completely. Problem 2 refers to the other cases that were not covered by Bennett in \cite{MB}.\\

{\bf Problem 1.}
{\it For both $a$ and $b$ odd, where $a>b>0$, find all solutions to the diophantine equation $3^{a}+3^{b}+2=(6j+2)^5$ or show that the only solution is $(a,b,j)=(3,1,0).$}\\

{\bf Problem 2.} 
{\it {For any positive integer $n$, with $(n,3)=1,$ find all solutions to $\gamma_{3}(n^{q})\leq3$ for $q$ a prime number where $q >1000$}}? \\

{\bf Problem 3.} 
{\it {For any positive integer $n$, with $(n,3)=1,$ we conjecture that $\gamma_{3}((6j+2)^5)>c\ln{(6j+2)^5}$ where $c$ is a constant such that $0<c<1/3\ln{(3)}.$}}\\

The binary analog for perfect powers has been well documented and can be found in \cite{MBM}.
\[
\]
{

\end{document}